\documentclass[a4paper]{article}
\usepackage[utf8]{inputenc} 
\usepackage[T1]{fontenc} 
\usepackage{lmodern} 
\usepackage{CJKutf8}
\usepackage{cite}
\usepackage[all]{xy}
\usepackage[dvips]{graphics,graphicx}
\usepackage{amsfonts,amssymb,amsmath,color,mathrsfs, amstext}
\usepackage{amsbsy, amsopn, amscd, amsxtra, amsthm,authblk}
\usepackage{enumerate,algorithmicx,algorithm,enumitem}
\usepackage{algpseudocode}
\usepackage[title]{appendix}
\usepackage{upref}
\usepackage{geometry,wrapfig}
\usepackage{amsthm,amsmath,amssymb}
\usepackage{relsize}
\usepackage{mathrsfs,mathtools}
\usepackage{bm,bbm}
\usepackage{ulem}
\usepackage{exscale}
\usepackage{tabularx}
\usepackage{threeparttable,multirow,dcolumn,booktabs}
\usepackage{algorithmicx,algorithm}
\usepackage{verbatim}
\usepackage{tabularx}
\usepackage{caption}
\usepackage{subfigure}
\usepackage[displaymath]{lineno}
\usepackage{float}
\usepackage{bm}
\usepackage{caption}
\usepackage{yhmath}
\usepackage{booktabs}
\usepackage{multirow}
\usepackage{makecell}
\usepackage[colorlinks,
            linkcolor=red,
            anchorcolor=red,
            citecolor=red
            ]{hyperref}
\usepackage{cleveref}

\usepackage[table,dvipsnames]{xcolor}

\definecolor{cgray}{gray}{0.80}
\newcolumntype{a}{>{\columncolor{cgray}}c}

\newcommand{\calD}{\mathcal{D}}

\newcommand{\TBD}[1]{\textcolor{black}{#1}}

\numberwithin{equation}{section}
\usepackage{epstopdf}
\usepackage{pdfpages}

\begin{document}

\title{EGPT-PINN: Entropy-enhanced Generative Pre-Trained Physics Informed Neural Networks for parameterized nonlinear conservation laws}

\author{
Yajie Ji\footnote{School of Mathematical Sciences, Shanghai Jiao Tong University, Shanghai 200240, China. Email: {\tt jiyajie595@sjtu.edu.cn}. Y. Ji acknowledges the support from the NSFC (No. 124B2023).},\quad
Yanlai Chen\footnote{Department of Mathematics, University of Massachusetts Dartmouth, North Dartmouth, MA 02747. Email: {\tt{yanlai.chen@umassd.edu}}. Y. Chen is partially supported by National Science Foundation grant DMS-2208277 and by Air Force Office of Scientific Research grant FA9550-23-1-0037.}, \quad  
Zhenli Xu\footnote{School of Mathematical Sciences, CMA-Shanghai and MOE-LSC, Shanghai Jiao Tong University, Shanghai 200240, China. Email: {\tt xuzl@sjtu.edu.cn}. Z. Xu acknowledges the support from the NSFC (grant Nos. 12325113 and 12426304) and the National Key R\&D Program of China (grant No. 2024YFA1012403).}
}

\date{}
\maketitle

\begin{abstract}
We propose an entropy-enhanced Generative Pre-Trained Physics-Informed Neural Network with a transform layer (EGPT-PINN) for solving parameterized nonlinear conservation laws. The EGPT-PINN extends the traditional physics-informed neural networks and its recently proposed generative pre-trained strategy for linear model reduction to nonlinear model reduction and shock-capturing domains. By utilizing an adaptive meta-network, a simultaneously trained transform layer, entropy enhancement strategies, implementable shock interaction analysis, and a separable training process, the EGPT-PINN efficiently captures complex parameter-dependent shock formations and interactions.  Numerical results of EGPT-PINN applied to the families of inviscid Burgers' equation and the Euler equations, parameterized by their initial conditions, demonstrate the robustness and accuracy of the proposed technique. It accurately solves the viscosity solution via very few neurons without leveraging any {\it a priori} knowledge of the equations or its initial condition. Moreover, via a simple augmentation of the loss function by model-data mismatch, we demonstrate the robustness of EGPT-PINN in solving inverse problems \TBD{more accurately than the vanilla and entropy-enhanced versions of PINN.}

\end{abstract}

\section{Introduction}
Hyperbolic conservation laws play a pivotal role in science and engineering thanks to their ability to describe the evolution of physical quantities such as mass, momentum, and energy \cite{dafermos1977generalized,book2005hyperbolic}. 
A primary challenge in numerically solving these equations is accurately ruling out nonphysical solutions and capturing discontinuities without introducing spurious oscillations. To address these challenges, numerous conventional numerical solvers, including finite difference \cite{LeVeque90,sod1978difference}, finite volume \cite{morton2007finite}, and discontinuous Galerkin methods \cite{cockburn2000development} have been developed over the past decades. 
Particularly significant developments include Godunov-type methods solving Riemann problems at cell interfaces, slope limiters combining high-order accuracy in smooth regions with non-oscillatory behavior near discontinuities, and essentially non-oscillatory and weighted essentially non-oscillatory  schemes \cite{shu2020WENO} adaptively selecting stencils based on the solution's smoothness. 

In recent years, machine learning algorithms have become valuable tools for solving complex Partial Differential Equations (PDEs) including those modeling fluid dynamics, electro-magnetics, heat transfer, and quantum mechanics. Leveraging the universal approximation capabilities of neural networks \cite{cybenko1989approximation, YAROTSKY2017103, Hanin2019}, researchers have introduced innovative methods like Physics-Informed Neural Networks (PINNs) \cite{raissi2019pinn} and similar approaches to efficiently tackle nonlinear conservation law equations \cite{zhang2022implicit, liu2024discontinuity, jagtap2020conservative, cai2022LSNN, mao2020physics, peyvan2024riemannonets}. However, the design of machine learning algorithms for hyperbolic conservation laws remains in an early stage of development and their model order reduction strategies essentially nonexistent. It continues to attract researchers' attention thanks to their promise of addressing some of the limitations found in traditional methods such as handling high-dimensional spaces \cite{EHanJentzen2017, HanJentzenE2018} and intricate mappings between parameters and solutions. 
Toward that end, operator learning approaches such as the Deep Operator Network (DeepONet) \cite{lu2021deeponet} and Fourier Neural Operator (FNO) \cite{liFNO}, have been introduced and extensively applied in approximating solutions or solution maps of PDEs and parametric PDEs (pPDEs). 
However, these methods overwhelmingly rely on a large amount of labeled data which is costly to generate in the PDE setting. The data and training do enable negligible evaluation time online given a new input. However, the direct evaluation of the trained network ignores the physical constraint enbodied by the PDE or pPDE.

The Reduced Basis Method (RBM)\cite{quarteroni2015reduced, HesthavenRozzaStamm2015, Haasdonk2017Review}  offers a promising philosophy by always embedding the physical constraints while significantly reducing the computational burden from having to repeatedly resolve pPDEs. The core idea behind the RBM is to rely on a mathematically rigorous greedy algorithm to construct a low-dimensional approximation space, from ground up, that captures the essential features of the solution space. This reduced space is generated using a set of carefully selected basis functions, derived from a series of precomputed high-fidelity solutions known as ``snapshots''. For any parameter values, the RBM then seeks a surrogate solution by typically satisfying the physical constraints in this reduced space. This allows for a rapid calculation when the so-called Kolmogorov $n$-width decays fast \cite{pinkus_n-widths_1985, BinevCohenDahmenDevorePetrovaWojtaszczyk, CohenDeVore2016, BuffaMadayPateraPrudhommeTurinici2011}. The RBM inspired the design of the GPT-PINN \cite{chen2024gpt}, the linear reduction regime in the PINN framework that features a tiny network with one sole hidden layer adaptively built with activation functions being fully pre-trained PINNs. 
However, when tackling transport-dominated problems, employing linear reduction methods such as GPT-PINN poses significant challenges. These challenges stem from the highly nonlinear and complex nature of transport phenomena, characterized by pronounced discontinuities that can evolve over time and exhibit sensitivity to parameter variations. 
In \cite{chen2024tgpt}, we present the Transformed Generative Pre-Trained Physics-Informed Neural Networks (TGPT-PINN), a novel framework employing nonlinear reduction strategies \cite{touze2021model,Welper2017}.  The method preserves the network structure and the unsupervised learning nature inherent in PINNs while providing a fast solver for pPDEs whose solution manifold features a slow-decaying Kolmogorov $n$-width \cite{OhlbergerRave2016, greif2019decay}.

While being able to accurately capture fully convective phenomena including that of full shocks with varying (linear or simple nonlinear) transport speed, the TGPT-PINN struggles in simulating inviscid fully nonlinear conservation laws. 
In this paper, we tackle that exact problem by proposing the Entropy-enhanced TGPT-PINN (EGPT-PINN) framework. 
These problems are challenging due to the emergence of dual challenges, the existence of non-physical solutions and the slow decay of the Kolmogorov $n$-width. 
In addition to an adaptive meta-network and a simultaneously trained transformation layer that are inherited from the TGPT-PINN, the EGPT-PINN introduces several key innovations to address this dual challenges for machine learning based model order reduction for parametric conservation laws. First, it replaces the original PDE with its characteristic form, simplifying the computation and enabling better capturing of shock waves. To address overfitting near discontinuities, the EGPT-PINN incorporates physics-dependent weights into the loss function, which helps the model focus on regions with shocks and improves solution accuracy. The loss function is further enhanced by including terms for the PDE, initial and boundary conditions, as well as the Rankine-Hugoniot (RH) condition \cite{jenny1998rankine,liu2024discontinuity}, which is essential for accurately modeling physically-relevant shock dynamics. Additionally, an indicator function is used to detect the location of discontinuities, allowing the model to adaptively concentrate on these critical areas and to predict shock intersection time. The EGPT-PINN similarly employs a nonlinear model order reduction approach by introducing a parameter-dependent transform layer, effectively handling parameter-dependent discontinuities and overcoming the limitations of linear model reduction in transport-dominated regimes. Moreover, the method adopts an offline-online structure. During the online stage, only the parameters of the transform layer and output layer are optimized significantly reducing computational effort. In the offline stage, a greedy algorithm is applied, generating a small number of neurons whose activation functions are  full PINNs pre-trained at judiciously selected locations in the parameter domain. These inherited features \TBD{from TGPT-PINN}, combined with the above-mentioned innovations, make EGPT-PINN a powerful and efficient tool for solving hyperbolic conservation laws with parameter-dependent discontinuities. To the best of our knowledge, this is the first attempt to accelerate the computing of parametric conservation laws leveraging nonlinear model reduction techniques in the neural network setting. To further illustrate the robustness of EGPT-PINN, we leverage it for inverse problems. The hallmark feature of PINNs and other neural network-based methods in solving inverse problems is that they do it with a computational cost that is comparable to the forward problem. The main tool is integrating the unknown parameters directly into the training process alongside network parameters. We demonstrate that this advantage translates from the full PINN to EGPT-PINN \TBD{when solving the challenging boundary inference problem for the Euler equation}.

The rest of this paper is organized as follows. Background materials including pPDEs, RBM, PINN, GPT-PINN and conservation laws are briefly presented in Section \ref{sec:bgd}. The main algorithm is given in Section \ref{sec:method} with numerical results for Burgers'  and Euler equations  in Section \ref{sec:numerics}. We draw conclusions in Section \ref{sec:conclusion}.

\section{Background}
\label{sec:bgd}

\subsection{pPDEs, PINN and GPT-PINN}
Parametric partial differential equations are extensively used in fields such as climate modeling, materials science, and biomedical engineering to describe complex systems. These parameters can represent physical properties, geometric configurations, or initial and boundary conditions. Solving these PDEs numerically is crucial for understanding and predicting the behavior of the complex systems they model. Recently, fast numerical algorithms for solving pPDEs have gained significant attention, particularly in engineering applications, which is driven by the need for repeated simulations of pPDEs in control, optimization, and design tasks. Consider the classic time-dependent PDE  on the spatial domain $\Omega \subset \mathbb{R}^d$ with boundary $\partial \Omega$ and parameter $\mu \in \mathcal{D}$:

\begin{subequations}
\label{eq:pPDE}
\begin{equation}
\label{eq:pPDE-a}
\frac{\partial}{\partial t} u(\bm{x}, t;\mu)+\mathcal{F}(u)(\bm{x}, t;\mu)=0, \quad \bm{x} \in \Omega, \quad t \in[0, T],
\end{equation}
\begin{equation}
\label{eq:pPDE-b}
\mathcal{G}(u)(\bm{x}, t;\mu)=0, \quad \bm{x} \in \partial \Omega, \quad t \in[0, T], 
\end{equation}
\begin{equation}
\label{eq:pPDE-c}
u(\bm{x}, 0;\mu)=u_0(\bm{x};\mu), \quad \bm{x} \in \Omega.
\end{equation}
\end{subequations}
Here $u(\bm{x},t;\mu)$ represents the solution, $u_0(\bm{x},\mu)$ provides the initial value, $\mathcal{F}$ is a differential operator and $\mathcal{G}$ encodes boundary information. 

Analyzing this system's behavior and its parameter dependence often requires executing thousands to millions of simulations of the PDE. Simulating it accurately and robustly is typically time-intensive, making the extensive repetition of simulations computationally infeasible with traditional numerical methods. To overcome this challenge, Model Order Reduction (MOR) has been gaining traction in the last few decades. 
The fundamental idea of the projection-based MOR \cite{BennerGugercinWillcox2015} is to construct a low-dimensional approximation of the solution space, and then find the Galerkin projection of the high-fidelity solution into this low-dimensional space \TBD{as a surrogate solution}. 

One notable example is the RBM which consists of two primary stages: the offline stage and the online stage. During the offline stage, high-fidelity solutions are generated for a carefully selected set of parameter values using methods like the proper orthogonal decomposition or the greedy algorithm. In the online stage, the reduced model can quickly solve for new parameter values using the pre-computed basis.

In recent years, PINNs \cite{raissi2019pinn} have demonstrated significant promise in addressing complex challenges in computational science and engineering. By embedding physical laws into the Deep Neural Networks (DNNs), PINNs offer a versatile and easy-to-implement approach to tackle a wide range of scientific and engineering problems. 
DNNs are composed of multiple layers of neurons that perform a series of linear transformations followed by nonlinear activations. This layered structure enables the practitioners to capture intricate patterns in data, making them particularly effective for tasks such as image recognition, natural language processing, and, more recently, solving differential equations. 
The output $\Psi_{\mathsf{NN}}({\bm x}; \Theta)$ of a DNN, parameterized by $\Theta$, can be represented as follows. 
\begin{equation}
\label{eq:DNN}
\Psi_{\mathsf{NN}}({\bm x}; \Theta) \coloneqq (C_L \circ \sigma \circ C_{L-1} \cdots \circ \sigma \circ C_1)(\bm x).
\end{equation}
It is then adopted as an approximation to the PDE solution $u(\bm x)\approx \Psi_{\mathsf{NN}}({\bm x}; \Theta)$. Here $\bm x$ is the input vector, $\Theta = \{\bm{W}_{\ell},\bm{b}_{\ell}\}_{{\ell}=1}^L$ with $\bm{W}_{\ell}$ and $\bm{b}_{\ell}$ respectively denoting the weight matrix and bias vector of the $\ell$-th ($1\leq \ell \leq L$) layer, and $\sigma$ denotes the activation. Each linear layer $C_{\ell}$ is defined by, 
\[
C_{\ell}(\bm x) = \bm{W}_{\ell} \bm x+\bm{b}_{\ell}.
\]

PINNs harness the expressivity of DNN to approximate solutions to PDEs by embedding the governing equation \eqref{eq:pPDE-a}, initial conditions \eqref{eq:pPDE-b}, and boundary conditions \eqref{eq:pPDE-c} into a loss function used for training the DNN. 
In the case of pPDE such as \eqref{eq:pPDE}, the resulting DNN is denoted as 
 $\Psi_{\mathsf{NN}}^{\mu}(\bm x,t; \Theta(\mu))$ to emphasize its parametric and time dependence. Whenever there is no confusion, we further simplify it as $\Psi_{\mathsf{NN}}^{\mu}(\bm x,t)$ or $\Psi_{\mathsf{NN}}^{\mu}$ for brevity. 
\begin{equation}
\label{PDE-loss}
\begin{aligned}
\mathcal{L}\left(\Psi_{\mathsf{NN}}^{\mu}(\bm x,t);\mu\right)&=\int_{\Omega \times (0,T]}\left\|\frac{\partial}{\partial t} \Psi_{\mathsf{NN}}^{\mu}(\bm x,t)+\mathcal{F}\left(\Psi_{\mathsf{NN}}^{\mu}(\bm x,t)\right)\right\|_2^2 d \bm{x}\, d t \\
&+ \int_{\Omega} \left\|\Psi_{\mathsf{NN}}^{\mu}(0,t)-u_0(\bm{x};\mu)\right\|_2^2 d \bm{x}+\int_{\partial \Omega \times [0,T]}\|\mathcal{G}\left(\Psi_{\mathsf{NN}}^{\mu}(\bm x,t)\right)\|_2^2 d \bm{x}\, d t.
\end{aligned}
\end{equation}
The training process minimizes this loss function through back propagation, incrementally adjusting the network's weights and biases, toward a better approximation with respect to the underlying physical constraints. PINNs represent a significant advancement in the numerical solution of PDEs thanks to the universal approximation capabilities of DNN \cite{cybenko1989approximation, YAROTSKY2017103, Hanin2019}, existence of software packages \cite{lu2021deepxde,haghighat2021sciann}, and the recent explosion of compute power \cite{agarwal2016tensorflow,wang2020survey}. However, when it comes to pPDEs, PINNs face the challenge of high computational cost associated with repeated training which is more pronounced than traditional methods due to the more resource-intensive nature of PINNs. 

To address this, researchers are exploring various strategies, such as transfer learning, where a trained network with one set of parameters is fine-tuned for others \cite{xu2023transfer}, and multi-fidelity methods that combine low-fidelity models with high-fidelity PINNs \cite{aliakbari2022multilowhigh,penwarden2022multifidelity}. In \cite{chen2024gpt}, the authors introduce a highly reduced neural network named GPT-PINN. This innovative approach reduces the size of PINNs required for unseen parameter values. The GPT-PINN is a network-of-networks, employing pre-trained PINNs as customized activation functions within the neurons of its single hidden layer. 
Using a mathematically rigorous greedy algorithm with residual-based error indicators, we select parameters $\mu^1,\mu^2,\cdots,\mu^n$ and obtain $\Psi_{\mathsf{NN}}^{\mu^1},\Psi_{\mathsf{NN}}^{\mu^2},\cdots,\Psi_{\mathsf{NN}}^{\mu^n}$ by training corresponding PINNs. For an untrained parameter $\mu$, we approximate the solution as:
\begin{equation}
u(\bm{x},t;\mu) \approx \Psi_{\mathrm{\mathsf{NN}}}^{\mu}(\bm x, t) \coloneqq \sum_{i=1}^n c_i(\mu) \Psi_{\mathrm{\mathsf{NN}}}^{\mu^i}(\bm x,t). 
\label{eq:gpt_ansatz_pinn}
\end{equation}
For problems featuring fast decay in Kolmogorov $n$-width, the method was shown to generate significant speedup.

\subsection{Conservative Hyperbolic PDEs}

Conservative hyperbolic PDEs describe conservation laws for quantities in physical systems and find wide applications across fluid mechanics, thermodynamics, and electromagnetics. They provide a mathematical framework for tracking the temporal and spatial evolution of these quantities. Typically, the governing equations read
\begin{equation}
\label{eq:Conlaw-a}
\frac{\partial\mathbf{U}(\bm{x},t)}{\partial t}+\nabla \cdot \mathbf{F(\bm{U})}=0, \quad \bm{x} \in \Omega, \quad t \in[0, T],
\end{equation}
where $\mathbf{U} = (u^1,u^2,\cdots,u^m)^T$ denotes the unknown function representing conserved quantities with initial condition $\phi(\bm{x})$ and an appropriate boundary condition. 
Additionally, $\mathbf{F} = (f^1,f^2,\cdots,f^m)^T$ is the flux function of the dimensions $m\times d$ of $\mathbf{U}$. 
 
In particular, the Riemann problem is fundamental in conservation laws, describing the evolution of discontinuous waves across physical states. It involves solving the initial value problem for one-dimensional conservation laws, where the initial condition is a discontinuous wave,
\begin{equation}
\label{eq:Riemann initial}
\phi(x) = \begin{cases}u_L, & x\leq x_c, \\ u_R, & x>x_c.\end{cases}
\end{equation}
The inherent complexity of dealing with discontinuities makes these problems particularly challenging to solve. 
However, their solutions are crucial for numerically resolving conservation laws since they form the algorithmic basis for many numerical methods.

The governing equation \eqref{eq:Conlaw-a} can be written in its characteristic form as
\begin{equation}
\label{eq:char_form}
\frac{\partial \mathbf{U}}{\partial t}+\lambda (\mathbf{U}) \frac{\partial \mathbf{U}}{\partial x}=0, \quad \bm{x} \in \Omega, \quad t \in[0, T],
\end{equation}
where $\lambda (\mathbf{U}) = \mathbf{F'(\bm{U})}$ is the Jacobian matrix of the flux function with respect to $\mathbf{U}$ and represents characteristic speed. 
This form transforms the system into equations along the characteristic curves 
which is defined by the ordinary differential equation
\[\frac{d\bm{x}}{dt} = \lambda (\mathbf{U}),\quad \bm{x}(0) = \bm{x}_0.\]
Along the characteristic curve, the solution can be formally written as:
\begin{equation}
\label{eq:implict_form}
    \mathbf{U}(\bm{x},t) = \mathbf{U}(\bm{x}-\lambda(\mathbf{U})t,0) = \phi(\bm{x}-\mathbf{F'(\bm{U})}t).
\end{equation}

It is easy to see that, if $\mathbf{F'(\bm{U})}$ is a constant, \eqref{eq:char_form} reduces to the transport equation. The solution manifold, as the constant speed changes, is of rank 1 after an appropriate shifting. This family of pPDEs can then be solved exactly in the model reduction setting by one neuron via a carefully designed algorithm identifying this shift without \textit{a priori} knowledge of (the constant) $\mathbf{F'(\bm{U})}$ \cite{chen2024tgpt}. 
On the other hand, if $\mathbf{F'(\bm{U})} = \mathbf{U}$, \eqref{eq:char_form} becomes the well-known Burgers' equation, a fundamental example of nonlinear wave phenomena involving discontinuous solutions. Additionally, the more complex Euler equations express the conservation of mass, momentum and energy, relating the velocity field $u$ and $v$ with the density field $\rho$ and the pressure field $p$\cite{book2005hyperbolic}. \TBD{Next, we will discuss each of the two equations in detail.}

{\it Burgers' equation} --
By choosing $\mathbf{U} = u$ and $\mathbf{F} =u^2 / 2$ in \eqref{eq:char_form}, the governing equation becomes 1D inviscid Burgers' problem,
\[
u_t+u u_x=0,\quad  x \in \Omega,\quad t\in [0,T].
\]

{\it Euler equations} --
By choosing $\mathbf{U} = (\rho,\rho u,E)^T$ and $\mathbf{F}=(\rho u,\rho u^2+p,u(E+p))^T$ in \eqref{eq:char_form} with $E = \rho u^2/2 + p/(\gamma-1)$ ($\gamma$ = 1.4 for ideal gas), the 1D Riemann problem for the Euler equation is formulated as:
\[
\frac{\partial \mathbf{U}}{\partial t}+\mathbf{A}\frac{\partial \mathbf{U}}{\partial x} = 0, \quad x \in \Omega,\quad  t\in [0,T].
\]
The matrix $\mathbf{A}$ is given by 
\[
\mathbf{A} = \mathbf{F}'(\mathbf{U}) =\left(\begin{array}{ccc}
0 & 1 & 0 \\
(\gamma-3)q & (3-\gamma) u & \gamma-1 \\
u\left(\frac{1}{2}(\gamma-1) u^2-H\right) & H-2(\gamma-1)q & \gamma u
\end{array}\right),
\]
where $q =u^2/2$ and $H =(E+p)/\rho.$
Thus, the governing equation in characteristic form is:
\begin{equation}
\left(\begin{array}{c}
\rho \\
\rho u \\
E
\end{array}\right)_t+\left(\begin{array}{ccc}
0 & 1 & 0 \\
(\gamma-3)q & (3-\gamma) u & \gamma-1 \\
u\left((\gamma-1)q-H\right) & H-2(\gamma-1)q & \gamma u
\end{array}\right)\left(\begin{array}{c}
\rho \\
\rho u \\
E
\end{array}\right)_x=0.
\end{equation}

\section{The EGPT-PINN algorithm}

In this section, we present the new Entropy-enhanced TGPT-PINN framework, designed to enhance the modeling of parameterized PDEs with complex discontinuities. By reformulating the PDE into its characteristic form and integrating a robust loss function that includes physics-based weights, the Rankine-Hugoniot condition, and an accurate prediction of shock intersection, the full model EGPT-PINN ensures higher fidelity in capturing shock dynamics and discontinuous behavior. A nonlinear model order reduction via a parameter-dependent transform layer addresses parameter-dependent discontinuities effectively, while an offline-online computational structure optimizes efficiency by separating parameterized training and fast inference stages. The following sections detail each of these advances and demonstrate how they contribute to both the accuracy and computational practicality of EGPT-PINN in solving transport-dominated conservation laws with shocks.

\label{sec:method}
\subsection{Physics-informed and entropy-aware full order model, \TBD{EPINN}}

The solutions of nonlinear conservation laws often develop discontinuities within finite time, even if the initial conditions are smooth \cite{liu2024discontinuity}. This typically results in a deterioration of solution accuracy near shocks and contact waves \cite{peyvan2024riemannonets}. The inherent complexity of handling discontinuities makes these problems particularly challenging to solve by traditional numerical methods \cite{LeVeque90}, warranting special care such as incorporation of artificial viscosity \cite{CockburnShuRKDG191, CockburnShuRKDG289, CockburnShuRKDG389} and adoption of slope limiting techniques \cite{CockburnShuReview, ZhangShu10max}. There is no exception for network-based approaches. 

\TBD{To make sure that the full PINN achieves reasonable accuracy, we build the entropy-aware loss function based on the weighted characteristic form of the equation and enforce the Rankine–Hugoniot condition. To differentiate from the vanilla version, we call the resulting entropy-enhanced PINN EPINN. The loss function of the EPINN drives the solution to the physically-relevant one during the optimization process, reflected by its $\mathcal{L}_{\rm{RH}}(u)$ term:}
\begin{equation}
\label{eq:PINN-loss}
\mathcal{L}(u)=\mathcal{L}_{\rm{int}}(u)+ \varepsilon_i \mathcal{L}_{\rm IC}(u)+\varepsilon_b \mathcal{L}_{\rm BC}(u)+\varepsilon_r \mathcal{L}_{\rm{RH}}(u).
\end{equation}
Here $\varepsilon_i$, $\varepsilon_b$, and $\varepsilon_r$ are parameters used to balance the four individual loss terms. For example, in shock tube problems, effectively training the initial conditions is crucial as they play a significant role in capturing the discontinuity within the domain. To achieve accurate solutions with rapid convergence, it becomes necessary to minimize the loss associated with the initial conditions at a faster rate compared to the loss of the weighted PDE. Therefore, we choose the weighting constants as $\varepsilon_i=10$ and $\varepsilon_b=10$. 

As regards to the four individual loss terms in \eqref{eq:PINN-loss}, $\mathcal{L}_{\rm IC}(u)$ and $\mathcal{L}_{\rm BC}(u)$ carry the standard form enforcing the boundary and initial conditions
\[
\mathcal{L}_{\rm IC}(u) \coloneqq 
\int_{\Omega} \left\|\mathbf{U}(\bm{x},0) - \phi(\bm{x})\right\|_2^2 d \bm{x},
\]
\[
\mathcal{L}_{\rm BC}(u) \coloneqq 
\int_{\partial \Omega \times [0,T]}\|  \mathbf{U}(\bm{x},t) - g(\bm{x},t)\|_2^2 d \bm{x}\, d t.
\]
The other two terms are explained in detail below, with the final component of this subsection devoted to shock interaction prediction.

\noindent {\bf Weighted characteristic form} --
To better capture shock waves and high-speed flows and to simplify computation, the original PDE is replaced with its characteristic form to describe the conservation laws \cite{zhang2022implicit,mao2020physics}. Furthermore, following \cite{liu2024discontinuity}, a physics-dependent weight is introduced into the loss function to mitigate overfitting that may arise from excessive training of the neural network near discontinuities. The PDE loss function for the proposed weighted characteristic form-based method is expressed as:
\begin{equation}
\mathcal{L}_{\rm int}(u) \coloneqq 
 \int_{\Omega \times (0,T]}\left\|\lambda(\bm{x},t)\,\left(\mathbf{U}_t+\mathbf{F'(\bm{U})} \mathbf{U}_{\bm{x}}\right)\right\|_2^2 d \bm{x}\, d t,
\end{equation}
where the factor $\lambda$, inspired by what was originally proposed by Liu {\it et al.} \cite{liu2024discontinuity}, is defined as:
\[
\lambda(\bm{x},t)=\frac{1}{\varepsilon_\lambda\,(|\nabla \cdot u|-\nabla \cdot u)+1}.
\]

\noindent {\bf Rankine–Hugoniot condition} -- Equation \eqref{eq:Conlaw-a} generally admits multiple weak solutions, necessitating additional conditions to select the physically correct solution, such as the Rankine-Hugoniot (RH) and other entropy conditions \cite{book2005hyperbolic}. 

The RH condition \eqref{eq:RH-cond} provides the relation between the shock speed $\mathbf{S}_{\rm RH}$ and the variables, as well as the flux across the discontinuities.
\begin{equation}
\label{eq:RH-cond}
    \mathbf{S}_{\rm RH} = \frac{\mathbf{F}(\mathbf{U_1})-\mathbf{F}(\mathbf{U_2})}{\mathbf{U}_1-\mathbf{U}_2},
\end{equation}
where $\mathbf{U}_1,\mathbf{F}(\mathbf{U_1}),\mathbf{U}_2,\mathbf{F}(\mathbf{U_2})$ are the conservative quantities and flux across the discontinuity respectively. We note that the shock speed for Burgers' equation can be further simplified by
\begin{equation*}
    \mathbf{S}_{{\rm RH},t} = \frac{\mathbf{F}(\mathbf{U}(\bm{x}+\Delta \bm x, t)-\mathbf{F}(\mathbf{U}(\bm{x}-\Delta \bm x, t))}{\mathbf{U}(\bm{x}+\Delta \bm x,t)-\mathbf{U}(\bm{x}-\Delta \bm x, t)} = \frac{\mathbf{U}(\bm{x}+\Delta \bm x, t)+\mathbf{U}(\bm{x}-\Delta \bm x, t)}{2}. 
\end{equation*}
The RH condition is only satisfied in the vicinity of discontinuities or strong shocks. Within a small time step $\Delta t$, the movement speed of discontinuities should adhere to the RH condition as follows:
\begin{equation}
\label{eq:RH-vec}
\bm{x}_{{\rm RH},t}+\mathbf{S}_{\rm RH}\cdot \Delta t = \bm{x}_{{\rm RH},t+\Delta t},   
\end{equation}
where $\bm{x}_{\rm RH,t}$ represents the location of discontinuity at time $t$. Therefore, we need an indicator to detect the location of discontinuities\cite{liu2024discontinuity}. 

For the Burgers' equation, the solution often shows discontinuities with steep gradients or significant differences between neighboring $\bm x$. We focus solely on the $u$ terms, as follows:
\begin{equation}
\label{eq:RH-B}
\lambda_{\rm RH,t}(\bm{x},\bm{x}-\Delta \bm x) = \begin{cases}
 \left |(u_1-u_2) \right |, & \text {if } \left | u_1-u_2 \right | >\varepsilon,\\
0,&\text{elsewhere}.
\end{cases}
\end{equation}
Once the indicator $\lambda_{\rm RH}$ is calculated, we are ready to define the penalty term for the RH condition for the Burgers' equation as:
\begin{equation}
\mathcal{L}_{\mathrm{RH}} \coloneqq 
 \int_{\mathcal{S}_{\rm RH}}\left\|\lambda_{{\rm RH}, t}\cdot\left(\bm{x}_{{\rm RH}, t}+\mathbf{S}_{\rm RH}\cdot {\Delta t} - \bm{x}_{{\rm RH},t+\Delta t}
 \right)\right\|_2^2 d \bm x \, dt.
\end{equation}

For the Euler equations, contact discontinuities are characterized by changes in density  without significant changes in pressure or velocity. In contrast, strong discontinuities involve substantial changes in both pressure, density and velocity. So for the RH-collocation points $(\bm{x},t)$ and adjacent points $(\bm{x}\pm \Delta \bm x,t)$, the indicator $\lambda_{\rm RH}$ functions as a filter to detect shock waves:
\begin{equation}
\label{eq:RH-ind}
\lambda_{{\rm RH},t}(\bm{x},\bm{x}-\Delta \bm x) = \begin{cases}
 \left | (p_1-p_2)(u_1-u_2) \right |, & \text {if } \left | p_1-p_2 \right | >\varepsilon_1 , \left | u_1-u_2 \right | >\varepsilon_2,\\
0,&\text{elsewhere}.
\end{cases}
\end{equation}
Here, $\varepsilon_1$ and $\varepsilon_2$ are two parameters used to detect jumps in shock waves, and their values can be adjusted depending on the specific problem at hand. Unless otherwise stated, we set $\varepsilon_1=\varepsilon_2=0.2.$ 
\TBD{The RH relation, coupling density velocity and pressure, can be simplified as \cite{liu2024discontinuity}:
\begin{equation}
\begin{aligned}
&\rho_1 \rho_2 \left(u_1-u_2\right)^2 =\left(p_1-p_2\right)\left(\rho_1-\rho_2\right), &\\
&\rho_1 \rho_2\left(e_1-e_2\right) =\frac{1}{2}\left(p_1+p_2\right)\left(\rho_1-\rho_2\right), &
\end{aligned}
\end{equation}}
where subscripts 1 and 2 denote the pre-shock and post-shock states respectively, dependent on $\lambda_{\rm RH}$ and computed by the adjacent points $(\bm{x}\pm \Delta \bm{x},t)$. 
Using the indicator $\lambda_{\rm RH}$ in \eqref{eq:RH-ind},  the penalty term for the RH condition in the Euler equations is defined as:
\begin{equation}
\begin{aligned}
\mathcal{L}_{\mathrm{RH}} \coloneqq 
&\left\|\lambda_{\rm RH}\left[\rho_1 \rho_2\left(u_1-u_2\right)^2 -\left(p_1-p_2\right)\left(\rho_1-\rho_2\right)\right]\right\|^2_2\\
+&\left\|\lambda_{\rm RH}\left[\rho_1 \rho_2\left(e_1-e_2\right)-\frac{1}{2}\left(p_1+p_2\right)\left(\rho_1-\rho_2\right)\right]\right\|^2_2.
\end{aligned}
\end{equation}

\noindent {\bf Shock interaction analysis and implementation} -- Simultaneously resolving multiple shocks by neural networks and capturing their merging during their propagation is a challenging task even with the two components above. The situation is further exacerbated by the need to track the subsequent propagation of the newly formed discontinuity. A vanilla PINN usually leads to an approximation error that is significantly larger after the merging than before. 
To address this, we start with the two-shock case and develop a novel PINN structure featuring {\it a shock merging-triggered restart mechanism}. This strategy involves two separate networks, each trained on distinct time regions while sharing the same spatial domain. The initial value for the second network is transferred from the final evaluation from the first network, with the first network's parameters not participating in the second network's training. 
To determine the time domain decomposition, we first apply a vanilla PINN within a smaller time interval $T_0$ to estimate the distance $L_{T_0}$ between the two shocks when reaching $T_0$, identifiable using the indicator in \eqref{eq:RH-B}. Given that the intersection of the shocks forms a triangle, we can leverage the properties of similar triangles to deduce the following relationship:
\begin{equation}
\frac{L_{T_0}}{L}=\frac{t_{\rm merge}-T_0}{t_{\rm merge}},
\label{eq:shockinter_predict}
\end{equation}
allowing us to solve for $t_{\rm merge}$.
Here $L$ is the initial separation between the shocks, and $t_{\rm merge}$ is the estimated merging time which is the location the time domain will be decomposed for the two networks.

\subsection{Reduced order model with a separable training process}

The TGPT-PINN approach \cite{chen2024tgpt} introduces an extension of the GPT-PINN in the context of nonlinear model reduction by employing a transform layer. This approach preserves the PINNs network structure and unsupervised learning nature while incorporating transformations for problems lacking (linear) low-rank structure. The TGPT-PINN retains the offline-online framework from GPT-PINN and the traditional RBMs, requiring the offline training of only a few PINNs for selected parameters to achieve high accuracy. The offline stage uses the greedy algorithm to determine these PINNs, while the online stage focuses on the rapid optimization of a small set of parameters. 

The EGPT-PINN inherits the structure and philosophy of TGPT-PINN while adopting the entropy-aware loss function \eqref{eq:PINN-loss} and a novel separable training. By introducing the transform layer in a GPT-PINN, the EGPT-PINN can capture the parameter-dependent discontinuity locations and achieves good accuracy. The transform layer $\mathcal{T}_{\mu,\eta}$ is designed as a mapping
$\mathcal{T}_{\mu,\mu^i}(\bm{x},t): \Omega \times [0, T] \longrightarrow \Omega \times [0, T],$ 
and we employ the linear transformation by Chen {\it et al.} \cite{chen2024tgpt} 
\begin{equation}
\mathcal{T}_{\mu, \eta}(\bm{x},t) \coloneqq 
{\rm Mod}_{\Omega, \mathcal{T}}
\left[
W_{\mu,\eta} 
\left(
\begin{tabular}{c}
     $\bm{x}$\\
     $t$ 
\end{tabular}
\right)
+ b_{\mu,\eta}
\right]
, \quad \eta = \mu^1, \dots, \mu^N.
\label{eq:transform}
\end{equation}
 Here, $W_{\mu,\eta} \in \mathbb{R}^{(d+1) \times (d+1)}$, and $b_{\mu,\eta} \in \mathbb{R}^{d+1}$, and ${\rm Mod}_{\Omega, \mathcal{T}}(\cdot)$ is an element-wise modulo map to ensure that each component of $\mathcal{T}$ outputs on the appropriate slice of $\Omega \times [0,T]$.

Following the procedure of the GPT-PINN, one uses the mathematically rigorous greedy algorithm to select parameters $\mu^1,\mu^2,\cdots,\mu^n$ and obtained $\Psi_{\mathsf{NN}}^{\mu^1}, \Psi_{\mathsf{NN}}^{\mu^2}, \cdots, \Psi_{\mathsf{NN}}^{\mu^n}$ by training corresponding PINNs in the offline stage. One can then approximate $\mathbf{U}(\bm{x},t;\mu)$ use the pre-trained PINNs as follows,
\begin{equation}
\mathbf{U}(\bm{x},t;\mu)\approx\Psi_{\mathrm{\mathsf{NN}}}^{\Theta(\mu)}(x, t) \coloneqq \sum_{i=1}^n c_i(\mu) \psi_{\mathrm{\mathsf{NN}}}^{\mu^i}(\mathcal{T}_{\mu,\mu^i}(x,t)),
\label{eq:tgpt_ansatz_pinn}
\end{equation}
where $\Theta(\mu) \coloneqq \left\{ \{W_{\mu, \mu^i}\}_{i=1}^n, \{b_{\mu, \mu^i}\}_{i=1}^n, \{c_i(\mu)\}_{i=1}^n \right\}$ represents the $n(d^2+3d+3)$ parameters to be trained in the EGPT-PINN, and  $n$ is the number of snapshots. 

\begin{figure}[h]
\centering
\includegraphics[width=\textwidth]{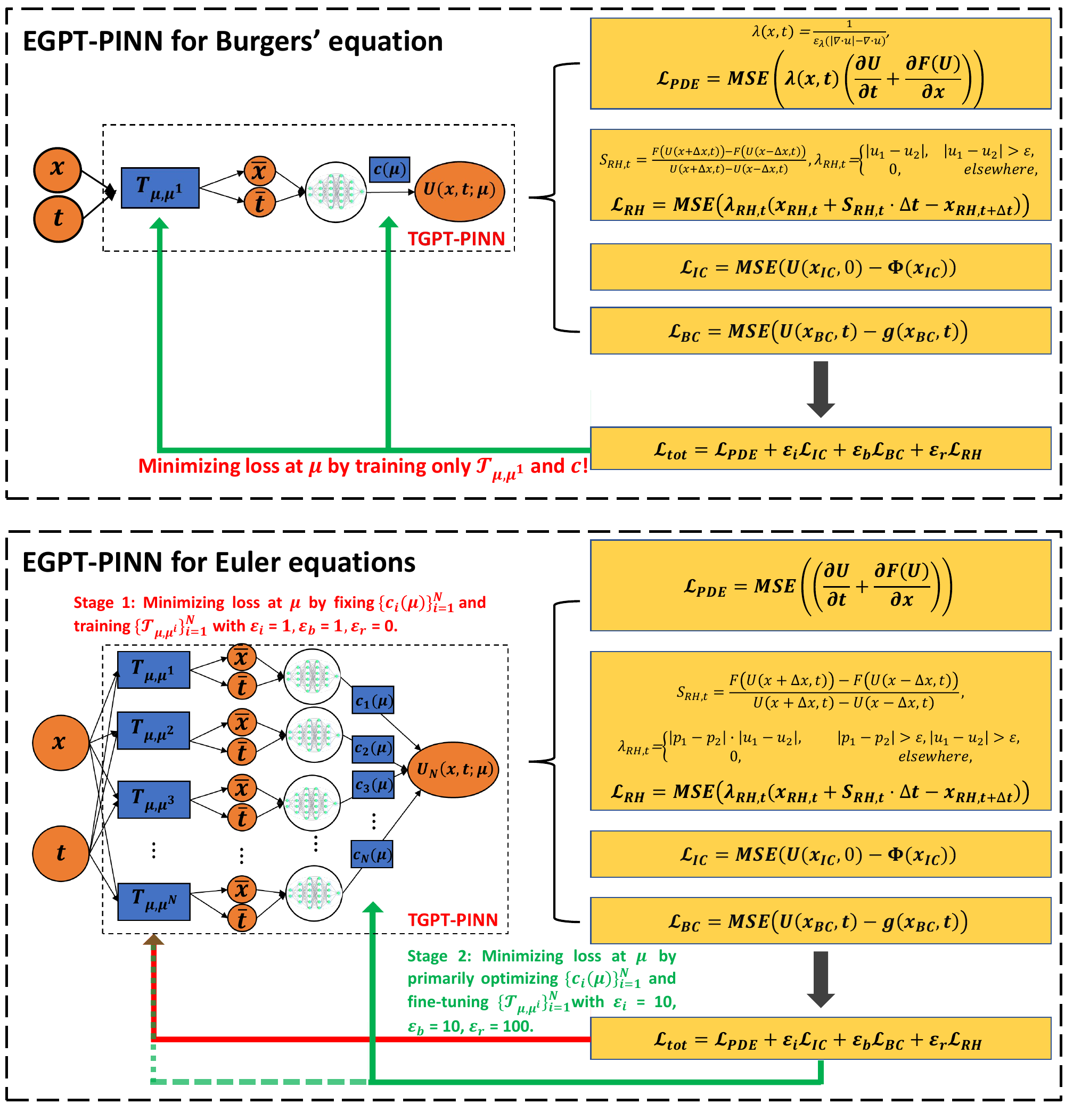}
  \caption{The EGPT-PINN design schematic for Burgers' (top) and Euler (bottom) equations. For any given parameter value $\mu$, a $\mu$-dependent loss is constructed and the coefficients $\{c_i(\mu)\}_{i=1}^n$ and the weights and biases in $\{\mathcal{T}_{\mu,\mu^i}\}_{i=1}^n$ are trained.}
\label{fig:EGPT-BurgersEuler}
\end{figure}

The loss function employed in the EGPT-PINN is defined in the same fashion as the full PINN, consisting of the governing equations, the initial value, the boundary condition and  the RH condition: 
\begin{subequations}
\label{eq:tgpt-pdeloss}
\begin{equation}
\label{tgpt-pdeloss-a}
\mathcal{L}_{\mathrm{int}}^{\mathrm{EGPT}}(\Theta(\mu))=\frac{1}{\left|\mathcal{C}_o^r\right|} \sum_{(\bm{x}, t) \in \mathcal{C}_o}\left\|\lambda(\bm{x},t)\,\left(\frac{\partial (\Psi_{\mathrm{\mathsf{NN}}}^{\Theta(\mu)})}{\partial t}+\frac{\partial (\mathbf{F}(\Psi_{\mathrm{\mathsf{NN}}}^{\Theta(\mu)}))}{\partial x} \Psi_{\mathrm{\mathsf{NN}}}^{\Theta(\mu)}\right)\right\|_2^2,
\end{equation}

\begin{equation}
\label{tgpt-pdeloss-b}
\mathcal{L}_{\mathrm{BC}}^{\mathrm{EGPT}} (\Theta(\mu))=\frac{1}{\mid \mathcal{C}_{\partial}^r \mid} \sum_{(\bm{x}, t) \in \mathcal{C}_{\partial}}\left\|\mathcal{G}\left( \Psi_{\mathrm{\mathsf{NN}}}^{\Theta(\mu)}\right)(\bm{x},t)\right\|_2^2,
\end{equation}

\begin{equation}
\label{tgpt-pdeloss-c}
\mathcal{L}_{\mathrm{IC}}^{\mathrm{EGPT}}(\Theta(\mu))=\frac{1}{\left|\mathcal{C}_i^r\right|} \sum_{\bm{x} \in \mathcal{C}_i}\left\| \Psi_{\mathrm{\mathsf{NN}}}^{\Theta(\mu)}(\bm{x}, 0)-u_0(\bm{x})\right\|_2^2,
\end{equation}

\begin{equation}
\label{tgpt-pdeloss-d}
\mathcal{L}_{\rm{RH}}^{\mathrm{EGPT}}(\Theta(\mu))=\frac{1}{\mid \mathcal{C}_{S}^r \mid} \sum_{(\bm{x}, t) \in \mathcal{C}_{S}}\left\|\lambda_{{\rm RH},t}\cdot\left(\bm{x}_{{\rm RH},t}+\mathbf{S}_{\rm RH}\cdot {\Delta t} - \bm{x}_{{\rm RH},t+\Delta t}
 \right)\right\|_2^2.
\end{equation}
\end{subequations}
The total loss function for the EGPT-PINN can be expressed by a combination of all these contributions:
\begin{equation}\label{eq:tgpt_loss}
\mathcal{L}_{\mathrm{PINN}}^{\mathrm{EGPT}}=\mathcal{L}_{\mathrm{int}}^{\mathrm{EGPT}} +\varepsilon_i\mathcal{L}_{\mathrm{IC}}^{\mathrm{EGPT}}+\varepsilon_b\mathcal{L}_{\mathrm{BC}}^{\mathrm{EGPT}}+\varepsilon_r \mathcal{L}_{\rm{RH}}^{\mathrm{EGPT}},
\end{equation}
where $\varepsilon_i$, $\varepsilon_b$ and $\varepsilon_r$ are balancing parameters for the loss terms, aligned with the values in \eqref{eq:PINN-loss}. By repeatedly updating  $\Theta(\mu)$ in the training process, the hyper-reduced network gradually converges towards the target solution. The algorithm details for the EGPT-PINN are presented in \Cref{fig:EGPT-BurgersEuler}.

The training process focuses on minimizing the loss function $\mathcal{L}_{\mathrm{PINN}}^{\mathrm{EGPT}}$ as defined in \eqref{eq:tgpt-pdeloss}. This is accomplished by utilizing standard techniques such as automatic differentiation and back propagation, employing the same learning rate for all the parameters in $\Theta(\mu)$. By iteratively updating the network's parameters based on calculated gradients, the EGPT-PINN learns to approximate the desired solution to the equations. 
However, the existence of discontinuities makes these problems challenging to optimize, especially for the Euler equations. 

We propose {\it a separable training technique} shown in \Cref{alg:gptpinn_twotraining}. The main purpose of this two-step process is to rely on the transform layer and give it sufficient time to better align shock waves. 
\begin{itemize}
\item At the first step, we freeze the output layer parameters $\{c_i(\mu)\}_{i=1}^n$ with $c_i(\mu) \equiv 1/n$ (or another initialization) and train the transform layer parameters $\{W_{\mu, \mu^i}\}_{i=1}^n$ and $\{b_{\mu, \mu^i}\}_{i=1}^n$ with a learning rate of $10^{-5}$. This step allows the network to identify and align shocks or discontinuities. The relatively larger updates to transform parameters allow the network focusing on capturing and accurately representing the locations and characteristics of these shocks, establishing a preliminary but reasonable approximation of the solution. 
\item In the second step of the training process, we set different learning rates for the output layer $\{c_i(\mu)\}_{i=1}^n$ and the transform layer $\{W_{\mu, \mu^i}\}_{i=1}^n$ and $\{b_{\mu, \mu^i}\}_{i=1}^n$, at $10^{-3}$ and $10^{-5}$ respectively.
\begin{itemize}
    \item The higher learning rate for the output layer facilitates rapid adjustments of the weights for different modes, allowing the network to converge to the desired solution more quickly. This helps in capturing the overall behavior and general trends of the problem. 
    \item The lower learning rate for the transform layer enables fine-tuning and precise adjustments. 
\end{itemize}
\end{itemize}
\begin{algorithm}[htbp]
    \caption{EGPT-PINN for parametric PDE: Online stage}
    \label{alg:gptpinn_twotraining}
    {\bf Input: }{A hyper-reduced initialization $\mathsf{NN}^n$ with ($\Psi_{\mathsf{NN}}^{\mu^1},\Psi_{\mathsf{NN}}^{\mu^2},\cdots,\Psi_{\mathsf{NN}}^{\mu^n}$), training set $\Xi_{\rm train} \subset \calD$}, learning rate $lr_1^{\text{in}},lr_2^{\text{in}},lr_1^{\text{out}}$ 
 and $lr_2^{\text{out}}$, max iterations $M_1$ and $M_2$, and tolerance $\delta$

\begin{algorithmic}[1]
\ForAll {$\mu \in \Xi_{\rm train}$ }
\While{\textit{iter}$\leq M_1$}
\State Train the $n$-neuron EGPT-PINN at $\mu$ with learning rate $lr = lr_1^{\text{in}}$ for the transform layer and $lr =lr_1^{\text{out}}$ for the output layer
\While{$\Delta_{\mathsf{NN}}^n(\Theta(\mu)) > \delta$ and \textit{iter}$+M_1 \leq M_2$}
\State Train the $n$-neuron EGPT-PINN at $\mu$ with learning rate $lr =lr_2^{\text{in}}$ for the transform layer and $lr =lr_2^{\text{out}}$ for the output layer
\State Record the indicator $\Delta_{\mathsf{NN}}^n(\Theta(\mu)) =\mathcal{L}_{\mathrm{PINN}}^{\mathrm{EGPT}}(\Theta(\mu)) $ in Eq. \eqref{eq:tgpt-pdeloss}

\EndWhile
\EndWhile
\EndFor
\end{algorithmic} 
 {\bf{Output:}~}  $n$-neuron EGPT-PINN $\mathsf{NN}^n(\mu)$
\end{algorithm}

Overall, this separable training allows the network to better capture the intricate details and complex patterns inherent to the PDE, leading to a more accurate and robust approximation of the solution. 

By integrating these different learning rates, the EGPT-PINN framework effectively balances the need for rapid alignment of the discontinuities and the requirement for accurate representation of the solution. This approach enhances the network's ability to capture both the global behavior and the local details, resulting in improved performance and higher accuracy in solving the problems.

Akin to the residual-based error estimation in traditional numerical solvers, our EGPT-PINN adopts the similar greedy offline process of the GPT-PINN \cite{chen2024gpt}. Indeed, we utilize the terminal loss $\mathcal{L}_{\mathrm{PINN}}^{\mathrm{EGPT}}$ for each parameter as the error indicator $\Delta_{\mathsf{NN}}^r(\Theta(\mu))$ which allows us to incrementally expand the EGPT-PINN hidden layer, in a ground-up fashion, from zero to (a pre-determined) $N$ neurons or until certain stopping criteria are met (e.g. error indicator falling below a threshold). At each step, the parameter value that is most poorly approximated by the current meta-network is selected. A full PINN is then pre-trained to augment the hidden layer. 
In this fashion, the meta-network learns the system's parametric dependencies one meta-neuron at a time. The offline training corresponds to the implementation of the greedy algorithm that is described in \Cref{alg:gptpinn_greedy}.

\begin{algorithm}[htbp]
    \caption{EGPT-PINN for parametric PDE: Offline stage}
    \label{alg:gptpinn_greedy}
    {\bf Input: }{A random (or given) $\mu^1$, training set $\Xi_{\rm train} \subset \calD$, and full EPINN}
\begin{algorithmic}[1]
\State Train a full EPINN at $\mu^1$ to obtain $\Psi_{\mathsf{NN}}^{\mu^1}$. Set $n=2$
\While{\textit{stopping criteria not met,}}
\State Train the ($n-1$)-neuron EGPT-PINN at $\mu$ for all $\mu \in \Xi_{\rm train}$ and record the indicator $\Delta_{\mathsf{NN}}^r(\Theta(\mu))$
 \State Choose $\mu^n = \displaystyle
  \mbox{\rm arg}\hspace*{-1pt}\max_{\mu\in{\Xi_{\rm train}}}
 {\Delta_{\mathsf{NN}}^r(\mu)}$
\State  Train a full PINN at $\mu^n$ to obtain $\Psi_{\mathsf{NN}}^{\mu^n}$. 
\State Update the EGPT-PINN by adding a neuron to the hidden EGPT-PINN layer to construct the $n$-neuron EGPT-PINN
\State Set $n \leftarrow n+1$
\EndWhile
\end{algorithmic} 
 {\bf{Output:}~}  $N$-neuron EGPT-PINN with {$N$ being the terminal index}
\end{algorithm}

\section{Numerical results}
\label{sec:numerics}

In this section, we perform numerical results to demonstrate the capability of the EGPT-PINN in effectively approximating the Riemann problems of the Burgers' and Euler equations parameterized by the initial conditions. Without relying on any prior knowledge or data, the method accurately captures, {\it often using only one neuron}, the complex dynamics and transitions of multiple shocks and rarefaction waves. 
The code for all these examples are published on GitHub at \href{https://github.com/DuktigYajie/VGPT-PINN}{https://github.com/DuktigYajie/EGPT-PINN}.

\subsection{Inviscid Burgers' equation}
\label{sec:Burgers}
For the 1D Burgers' equation, our extensive tests include five cases as listed in Table \ref{tab:BurgersTests}. They are one single shock propagating by itself, two shocks becoming one, a smooth initial condition developing into a shock, one rarefaction wave, and finally a shock interacting with a rarefaction wave. They are named ${\rm B_{1S}}$,  ${\rm B_{2S}}$, ${\rm B_{Sm}}$, ${\rm B_{R}}$ and ${\rm B_{RS}}$ respectively. \TBD{These five test cases represent basic formation and interaction examples of shocks and rarefaction waves, and are building blocks for more complicated scenarios.} 

\begin{table}[htb]
\begin{tabular}{|c|l|l|}
\hline
Test name & Initial condition $\phi(x, t = 0)$ & Boundary condition $\mathcal{G}(\partial \Omega, t)$ \\
\hline
${\rm B_{1S}}$ & 
$\phi(x, 0) = \begin{cases}\mu, & x \leqslant 0 \\ 
0, & x>0\end{cases}$
&
$\begin{cases}
 g(-1,t) = \mu,\quad t\in[0,1]\\
 g(1,t)  = 0,\quad t\in [0,1]
\end{cases}$\\
\hline
${\rm B_{2S}}$ & 
$\phi(x,0)= \begin{cases}\mu_2, & x \leqslant-1 / 2, \\ \mu_1, & -1 / 2<x \leqslant 1 / 2 \\ 0, & x>1 / 2\end{cases}$
&
$\begin{cases}g(x=-1, t)=\mu_2,\quad t \in(0,3 / 2]\\
g(x=3 / 2, t)=0,\quad t \in(0,3 / 2]
\end{cases}$\\
\hline
${\rm B_{Sm}}$ & 
$\phi(x,0)=\mu_1\sin(2\pi x) + \mu_2$
&
$g(0,t)=g(1,t),\quad t\in(0,1]$
\\
\hline
${\rm B_{R}}$ &
$\phi(x,0)= \begin{cases}0, & x \leqslant 0 \\ 
\mu, & x>0\end{cases}$
&
$\begin{cases}
 g(-1,t) = 0,\quad t\in [0,1]\\
 g(1,t)  = \mu,\quad t\in [0,1]
\end{cases}$\\
\hline
${\rm B_{RS}}$ & 
$\phi(x,0)= \begin{cases}0, & x \leqslant-1 / 2 \\ 
\mu_1, & -1 / 2<x \leqslant 1 / 2 \\ 
0, & x>1 / 2\end{cases}$
&
$u(-1, t)=u(3 / 2, t)=0$\\
\hline
\end{tabular}
\caption{Five test cases for the inviscid Burgers' equation}
\label{tab:BurgersTests}
\end{table}

The underlying EPINNs have 5 hidden layers with 20 neurons each for all tests. They are trained with the Adam optimizer, using 25,000 epochs and an initial learning rate of 0.001. The space-time collocation set  is a  $99 \times 199$ uniform grid excluding the boundaries. Moreover, 100 points each are placed to calculate the loss from the boundary and initial conditions. \TBD{The EGPT-PINN is optimized using the Adam optimizer with a learning rate of 0.001 and trained over 1,000 epochs, which is less than one-tenth of the epochs required for the full EPINN training.}

\begin{figure}[htbp]
\includegraphics[width=\textwidth]{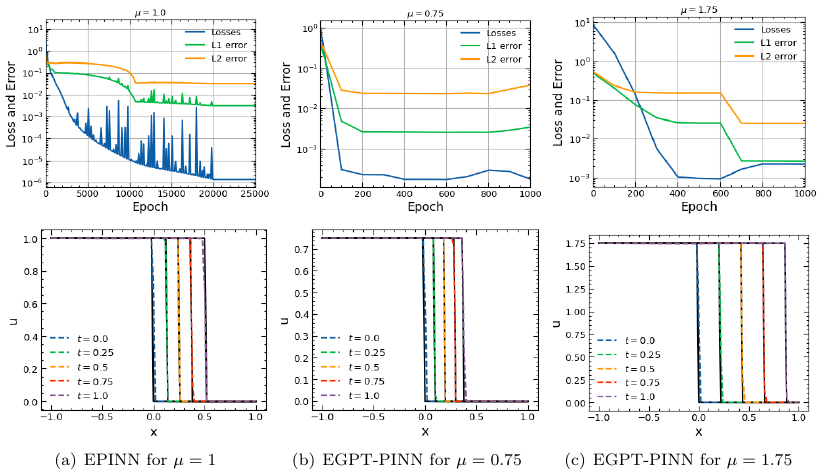}
  
\caption{Results of Case ${\rm B_{1S}}$. (a) Full EPINN loss at $\mu = 1.0$ and comparison with the exact solution; (b) One-neuron EGPT-PINN loss and comparison with the exact solution at $\mu = 0.75$; (c) One-neuron EGPT-PINN loss and comparison with the exact solution at $\mu = 1.75$.}
\label{fig:shock1}
\end{figure}

{\bf Case ${\rm B_{1S}}$} -- We choose the spatial domain $\Omega=[-1,1]$, the time interval $t \in [0,1]$, and the parameter domain $\mathcal{D}=[1,2]$. The weak solution of the Burgers' equation exhibits discontinuity along the straight line $x=\mu t/2$:
\[
u(x, t)= \begin{cases}\mu, & x<\mu t / 2, \\ 0, & x>\mu t / 2.\end{cases}
\]
We pre-train a single EPINN at $\mu = 1$. The full-EPINN solution, its comparison with the exact solution, and the reduction of the loss function are presented in \Cref{fig:shock1} (a) \TBD{attesting an accuracy of two significant digits}. 

\begin{figure}[ht]
\centering
\includegraphics[width=\textwidth]{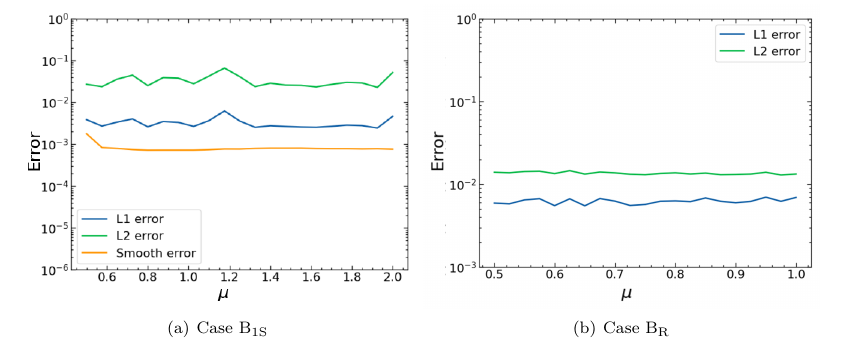}

\caption{Error test across the parameter domain for case ${\rm B_{1S}}$ (a) and ${\rm B_{R}}$ (b).}
\label{fig:shock1+4}
\end{figure}

The one-neuron EGPT-PINN is tested on $\mu \in [1:0.025:2]$. We evaluate the L1 and L2 errors across the parameter domain and present the results in \Cref{fig:shock1+4} (a). The variation of solutions, detailed errors and losses for two particular parameters unseen during training are shown in \Cref{fig:shock1} (b,c). Comparing with the full-EPINN results, it is clear that our EGPT-PINN solutions achieve the same accuracy. Note that here the ``smooth error'', \TBD{one order of magnitude lower than the whole domain error,} is calculated by excluding space-time points within a distance of $0.02$ from the discontinuity line.

\begin{figure}[htbp]
  \centering
\includegraphics[width=\textwidth]{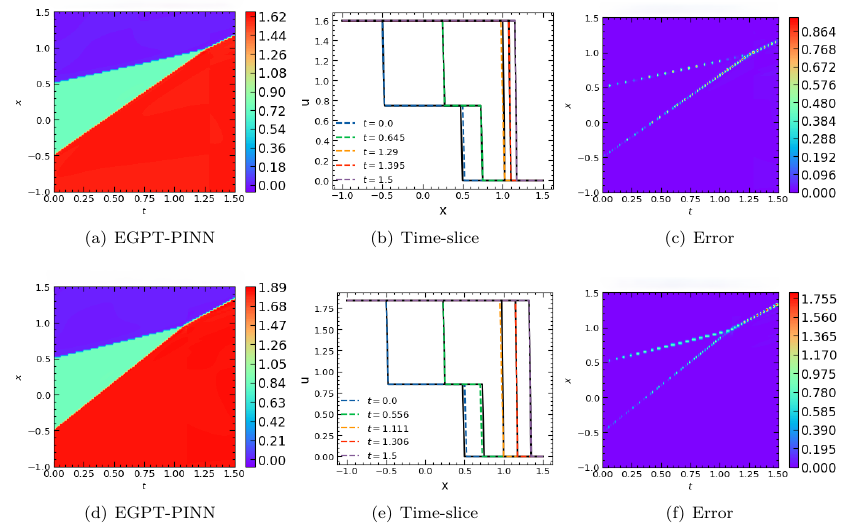}
  \caption{Results of Case ${\rm B_{2S}}$. (a-c) Two-neurons EGPT-PINN results for $\mu = (0.75,1.6)$; (d-f) Two-neurons EGPT-PINN results for $\mu =(0.85, 1.84)$.}
  \label{fig:shock2}
\end{figure}

{\bf Case ${\rm B_{2S}}$} -- Here, we set $\Omega=[-1,3 / 2]$ and 
parameter domain $\mathcal{D} = [0.5, 1.0] \times [1.6, 2.0]$. By using the RH jump condition, there exist two discontinuous lines $x=((\mu_1+\mu_2) t-1) / 2$ and $x=(\mu_1 t+1) / 2$ up to time $t=2/\mu_2$ when these two shocks  merge into one. The correct viscosity solution in this case is
\[
u(x, t)= \begin{cases}\mu_2, & t<2/\mu_2, x<((\mu_1+\mu_2)t-1) / 2, \\ \mu_1, & t<2/\mu_2,((\mu_1+\mu_2)t-1) / 2<x<(\mu_1t+1) / 2, \\ 0, & t<2/\mu_2, x>(\mu_1t+1) / 2, \\ \mu_2, & t>2/\mu_2, x<(\mu_2/2)t+(\mu_1/\mu_2-1/2), \\ 0, & t>2/\mu_2, x>(\mu_2/2)t+(\mu_1/\mu_2-1/2).\end{cases}
\]
We generate two independently-trained networks that are connected at the meeting time whose exact value $t=2/\mu_2$ is numerically predicted by the networks via \eqref{eq:shockinter_predict}. The first-stage is trained using a learning rate of 0.001 and 60,000 epochs, while the second-stage has a learning rate 0.0001 and 40,000 epochs. The EGPT-PINN uses the same temporal segmentation network structure as the EPINN, dividing the time domain into segments with separate networks trained for each segment. This structure facilitates effective handling of time-dependent problems through division of the time domain into manageable segments, each with its own trained network. The EGPT-PINN matches EPINN's learning rates but requires less than one-tenth of the epochs. The results of EGPT-PINN for two unseen parameters are shown in \Cref{fig:shock2}. It is clear that the shocks and their interaction are accurately captured by our new method. Moreover, the EGPT-PINN errors match those of the PINN's \TBD{with both achieving two digits of accuracy.}

{\bf Case ${\rm B_{Sm}}$} -- In fact, even smooth and continuous initial conditions can lead to shock formation. 
We consider the 1D inviscid Burgers' equation with space domain $\Omega=[0,1]$, time domain $t\in[0,1]$ and parameter domain $\mathcal{D} = [0.5,1.5] \times [0.1,0.4]$, with a smooth initial value and periodic boundary condition. 
In this scenario, the exact solution cannot be explicitly expressed. To obtain a reference solution, we adopt the Godunov scheme \cite{LeVeque90} and solve the equation on a sufficiently fine mesh. 

\begin{figure}[htbp]
\centering
\includegraphics[width=\textwidth]{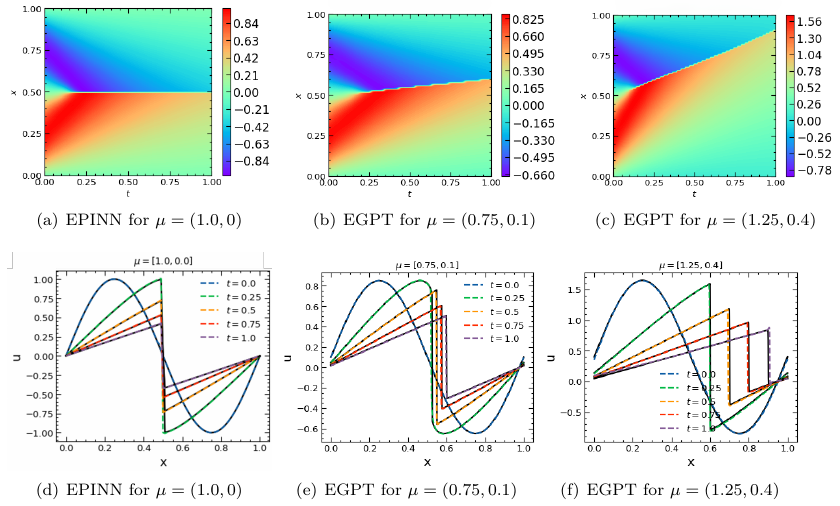}
\caption{Results of Case ${\rm B_{Sm}}$. (a, d) Full EPINN results for $\mu = (1.0, 0)$ , and the one-neuron EGPT-PINN results for $\mu = (0.75, 0.1)$ (b, e), $\mu = (1.25, 0.4)$ (c, f).}
\label{fig:shock3}
\end{figure}

The full-EPINN solution and this reference solution are shown in \Cref{fig:shock3} (a,d). The EGPT-PINN solutions and errors under three unseen parameters are shown in \Cref{fig:shock3} (b,e) and (c,f). \TBD{It is clear that the EGPT-PINN scheme captures the location of the shock very accurately, even when it is moving, by only one neuron. Moreover, we achieve better than 2-digit accuracy away from the shock.}

\begin{figure}[htbp]
\centering
\includegraphics[width=\textwidth]{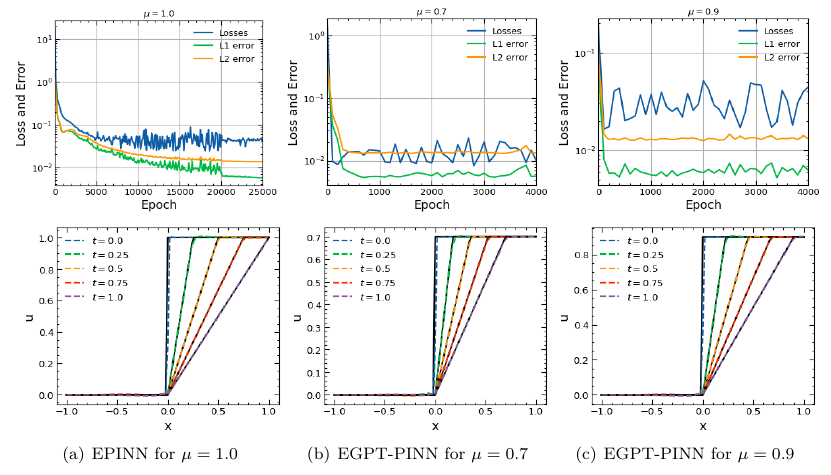}

\caption{Results of Case ${\rm B_{R}}$. Full EPINN results for $\mu = 1.0$ (a), one-neuron EGPT-PINN results for two unseen parameters $\mu = 0.7$ (b) and $\mu = 0.9$ (c).}
\label{fig:rare}
\end{figure}

{\bf Case ${\rm B_{R}}$} -- This is the rarefaction-wave case.
We take $\Omega=[-1,1]$. To determine a physically relevant solution, entropy condition, 
$
f^{\prime}\left(u^{-}\right)>\xi^{\prime}(t)>f^{\prime}\left(u^{+}\right),
$
must be enforced.  This leads to the exact solution
$$
u(x, t)= \begin{cases}0, & \text { for } x<0, \\ x / t, & \text { for } 0 \leqslant x \leqslant t, \\ \mu, & \text { for } x>t.\end{cases}
$$
The full EPINN solutions, errors, and the reduction of the loss function when $\mu=1$ are shown in \Cref{fig:rare}(a). Subsequently, we test the EGPT-PINN formed with this single neuron for $\mu \in [0.5:0.0125:1]$ and show the L1 error and L2 error in Figure \ref{fig:shock1+4} (b). The variation of solutions, errors and losses for two random parameters are shown in Panels (b, c) of this figure. \TBD{Once again, the single-neuron EGPT-PINN solutions achieve an accuracy comparable to that of the full PINN. The rarefaction of what is initially a shock is precisely captured as time evolves. We note that, had the entropy condition not been enforced, the accuracy of the solution would have deteriorated by multiple orders of magnitude. This underscores the robustness of the scheme due to its incorporation of the entropy condition.}

\begin{figure}[htbp]
\centering
\includegraphics[width=\textwidth]{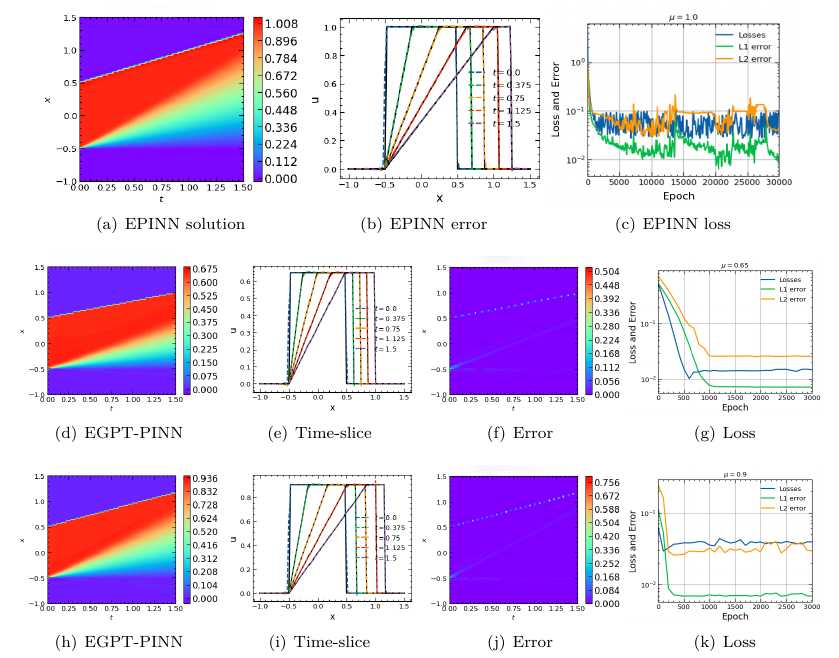}
\caption{Results of Case ${\rm B_{RS}}$. Full EPINN results for $\mu = 1.0$ (a-c). One-neuron EGPT-PINN results for $\mu = 0.5$ (d-g),  $\mu = 0.9$ (h-k).}
\label{fig:shockrare}
\end{figure}

{\bf Case ${\rm B_{RS}}$} -- The fifth is the case of the rarefaction-shock interaction. 
We choose $\Omega=[-1,3 / 2]$ and parameter domain $\mathcal{D} = [0.5,1.0]$. 
The RH jump condition tells us that there exists a discontinuity at $x=(\mu_1 t+1) / 2$ and a rarefaction region. Therefore, the viscosity solution in this case satisfying the entropy condition reads
$$
u(x, t)= \begin{cases}0, & x\geq (\mu_1t+1) / 2, \\ \mu_1, &\mu_1t-1/ 2\leq x<(\mu_1t+1)/ 2, \\ x/t, & -1/2\leq x<\mu_1 t-1/ 2, \\ 0, & x\leq-1/2. \end{cases}
$$
The PINN solutions, errors, and the reduction of the loss function for equation at $\mu=1$ are shown in \Cref{fig:shockrare} (a-c). The variation in EGPT-PINN's solutions, errors and loss under two unseen parameters are shown in \Cref{fig:shockrare} (d-k). \TBD{These results demonstrate that the EGPT-PINN with one neuron is capable of accurately capturing two different kinds of waves that are simultaneously evolving.}

\newpage
\subsection{Euler equations}
\label{sec:Euler-1D}
For the 1D Euler equations, we build the underlying EPINNs featuring 6 hidden layers with 60 neurons each for two tests, the classical Sod and Lax problems. 
The EGPT-PINN with $N$-snapshots has single hidden layer with $N$ neurons, using the two step training in  \Cref{alg:gptpinn_twotraining} with Adam optimizer and an initial learning rate of 0.001.

\begin{figure}[H]
\centering
\includegraphics[width=\textwidth]{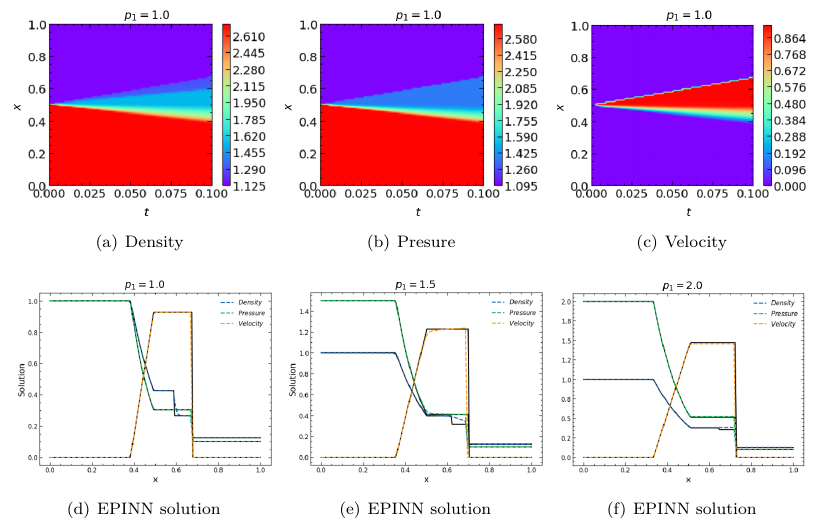}
  \caption{EPINN results for the Sod problem. (a-c) $x-t$ solution for $p_1 = 1.0$, (d-f) EPINN solution at final time and its comparison with a reference solution.}
  \label{fig:Full-PINN-Sod}
\end{figure}

{\bf Sod problem} -- 
The Sod problem\cite{sod1978difference} is a one-dimensional Riemann problem characterized by initial constant states within a tube of unit length. The parametric initial condition for  Eq.~\eqref{eq:char_form} is given by
\begin{equation}
    (\rho, u, p)= \begin{cases}(1,0,p_1), & \text {for} \quad 0 \leq x \leq 0.5, \\ (0.125,0,0.1), & \text { for } \quad 0.5<x \leq 1,\end{cases}
\end{equation}
where $p_1\in [1.0,2.0].$ 
The collocation set consists of randomly selected 5,000 interior points. We place 100 points each for the initial and boundary conditions. Additionally, 100 RH collocation points are drawn from a uniform mesh of $100 \times 200$ in the $X \times T$ space. The stopping criterion for the loss was set at $10^{-5}$ with a maximum of 30,000 epochs. In the underlying EPINN, we set
$\varepsilon_i=\varepsilon_b=10$ and $\varepsilon_r=100$
in the loss function \eqref{eq:PINN-loss}. The three components of the space-time EPINN solution for $p_1 = 1.0$ are shown in   \Cref{fig:Full-PINN-Sod} (a-c). The EPINN and exact solutions at final time $T=0.1$ for three different parameters are shown in  \Cref{fig:Full-PINN-Sod} (d-f).

\begin{figure}[!h]
  \centering
\includegraphics[width=\textwidth]{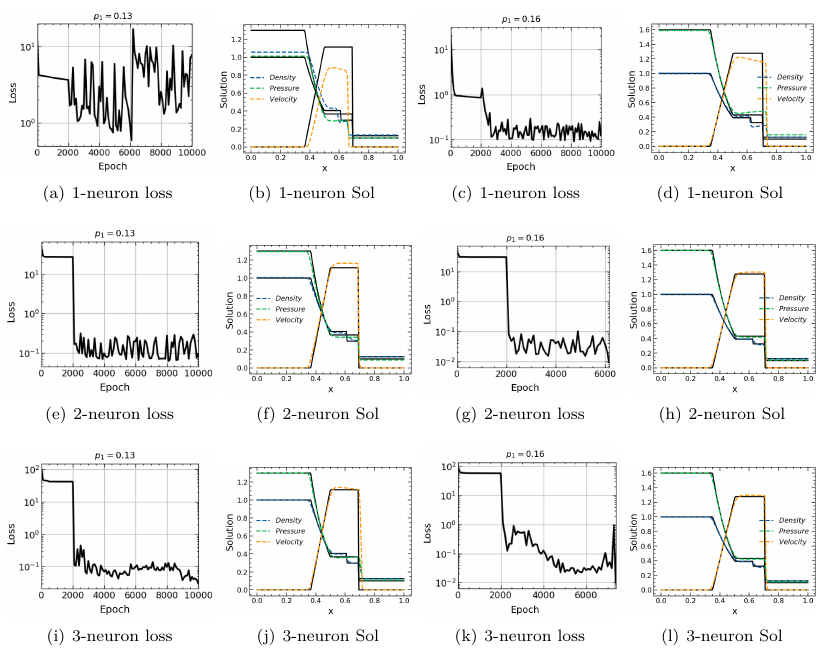}
  \caption{EGPT-PINN results for the Sod problem. Loss reduction and solution comparison for $\mu = 1.3$ (left two columns) or $\mu = 1.6$ (right two columns) for the EGPT-PINN generated by $1, 2$ and $3$ neurons (from top to bottom).}
  \label{fig:TGPT-PINN-Sod}
\end{figure}

In the EGPT-PINN process, we only need 3,000 interior points from a uniform mesh of $100 \times 200$ in the $x \times t$ space and utilized 100 initial points, 100 boundary points, and 100 RH collocation points. We set $\lambda(x,t) \equiv 1,\varepsilon_i=\varepsilon_b=1$ and $\varepsilon_r=0$ in the loss function of the first stage, and $\varepsilon_i=\varepsilon_b=10$ and $\varepsilon_r=100$ in the second stage. We consider all Sod problems with parameter domain $[1.0, 2.0] \ni p_1$ and test EGPT-PINNs with 1 to 3 neurons. As demonstrated by \Cref{fig:TGPT-PINN-Sod}, the approximate solution at $t = 0.1$ gradually approaches the reference solution. Using just three neurons, the approximate solution attains high accuracy, capturing the shock waves precisely. The EGPT-PINN solutions and the loss are shown in \Cref{fig:TGPT-PINN-Sod}. \TBD{It is clear by comparing with \Cref{fig:Full-PINN-Sod} that, with only three neurons, the EGPT-PINN solutions achieve an accuracy that is indistinguishable from that of the full EPINNs.}

\begin{figure}[thbp]
\centering
\includegraphics[width=\textwidth]{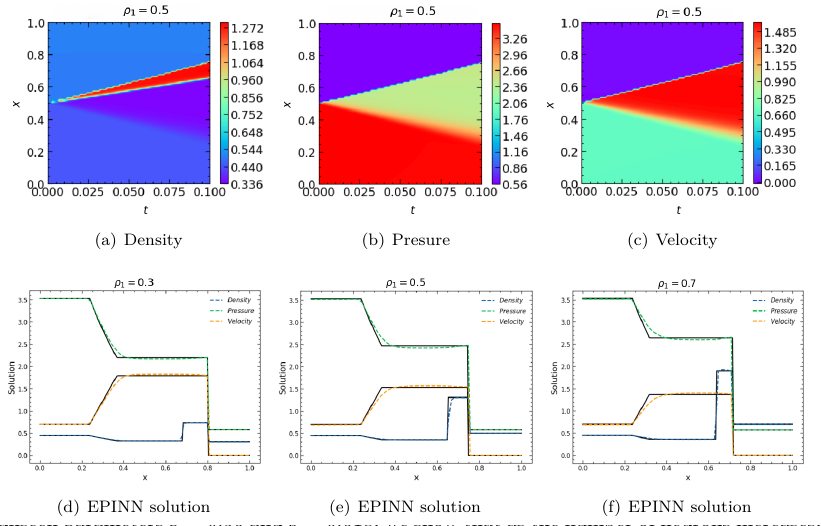}
  \caption{EPINN results for the Lax problem. (a-c) $x-t$ solution for $\rho_1 = 0.5$, (d-f) EPINN solution at final time and its comparison with a reference solution.}
  \label{fig:Lax}
\end{figure}

{\bf Lax problem} --
The Lax problem\cite{lax1957asymptotic} is another Riemann problem that contains a strong shock and strong contact.  These features create more complex wave interactions and require more sophisticated numerical techniques to resolve accurately. We endow the following parametric initial condition
\begin{equation}
    (\rho, u, p)= \begin{cases}(0.445,0.698,3.528), & \text { if } 0 \leq x \leq 0.5, \\ (\rho_1,0,0.571), & \text { if } 0.5<x \leq 1 .\end{cases}
\end{equation}
where $\rho_1\in [0.3,0.7].$ 
We set $\varepsilon_i=\varepsilon_b=10$ and $\varepsilon_r =100$ in the loss function. The collocation set consists of randomly selected 30,000 interior points. Moreover, 1000 initial points, 1000 boundary points and 1000 RH collocations with Latin Hypercube Sampling  method \cite{loh1996latin} in the $x \times t$ space, and the total loss was 0.01 after training for 20,000 epochs. The EPINN solutions for Lax problem with $\rho = 0.5$ are shown in \Cref{fig:Lax} (a-c). We also present the well-trained full-EPINNs for some parameters and compare them with the reference solution in \Cref{fig:Lax} (d-f). 

\begin{figure}[htbp]
\centering
\includegraphics[width=\textwidth]{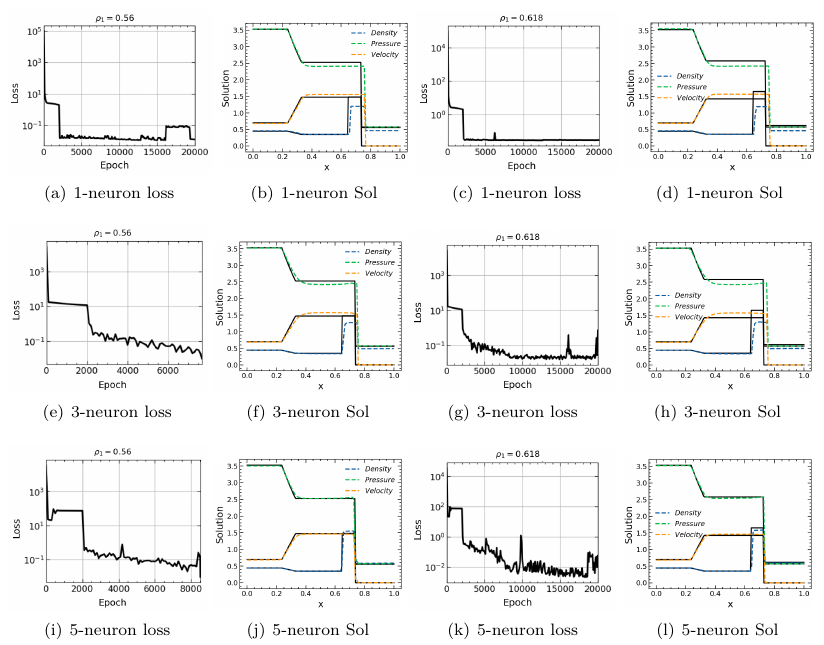}  
  \caption{EGPT-PINN results for the Lax problem: Loss reduction and solution comparison for $\mu = 0.56$ (left two columns) or $\mu = 0.618$ (right two columns) for the EGPT-PINN generated by $1, 3$ and $5$ neurons (from top to bottom).}
  \label{fig:TGPT-PINN-Lax}
\end{figure}

In the EGPT-PINN process, we set $\lambda(x,t) \equiv 1, \varepsilon_i=\varepsilon_b=1$ and $\varepsilon_r=0$
for the loss function in the first stage, and then
$\varepsilon_i=\varepsilon_b=1$ and $\varepsilon_r=10$ in the second stage. Using the unseen parameters $\rho = 0.56$ and $\rho = 0.618$, we show that as the number of neurons increases, the EGPT-PINN solution steadily approaches the reference solution in \Cref{fig:TGPT-PINN-Lax}. In the end, using just five neurons, the approximate solution attains high accuracy, capturing both the shock waves and shock contact precisely.

\subsection{Inverse problem}
To illustrate the robustness of EGPT-PINN, we task it to solve inverse problems. Toward that end, we consider the general formulation of an inverse problem: Given partial information about the solution, the objective is to infer unknown parameters or initial conditions of the underlying system. It has been widely demonstrated that EPINNs possess unique advantages in solving inverse problems. Traditional approaches often require multiple iterative steps, leading to significantly higher computational costs compared to forward problems. In contrast, neural network-based methods integrate the unknown parameters directly into the training process alongside network parameters, thereby achieving computational costs comparable to solving forward problems. The question we strive to answer in this paper is whether this advantage translates from the full EPINN to EGPT-PINN.

Taking the Sod shock tube problem as an example, suppose that we obtain measurements of density $\rho^*$, pressure $p^*$, and velocity $u^*$ at a specific location $(x^*,t^*)$ through experimental observations. The goal is to determine the corresponding initial conditions, such as the initial pressure. By incorporating the observed solution $\mathbf{U}^*=(\rho^*,u^*,p^*)$ into the loss function and treating the initial conditions as trainable parameters, the network can simultaneously optimize the unknown parameters and fit the desired initial conditions during training, ultimately yielding an optimal solution. We achieve this by simply augmenting the forward loss function with a data-solution mismatch
\begin{equation}
\mathcal{L}_{\mathrm{PINNInv}}^{\mathrm{EGPT}}=\mathcal{L}_{\mathrm{PINN}}^{\mathrm{EGPT}}+\varepsilon_d \left\|\mathbf{U}(x^*,t^*)-(\rho^*,u^*,p^*) \right\|_2,
\end{equation}
where $\mathcal{L}_{\mathrm{PINN}}^{\mathrm{EGPT}}$ the standard EGPT-PINN loss introduced in \eqref{eq:tgpt_loss}. In this setting, we consider the same configuration as in \Cref{sec:Euler-1D} for the Sod problem but assume that the initial pressure $p_1$ is unknown. Using EGPT-PINN, we simultaneously train both the parameterized solution $\Theta(\mu)$ and the unknown initial pressure $p_1$.

\begin{figure}[htbp]
    \centering
\includegraphics[width=\textwidth]{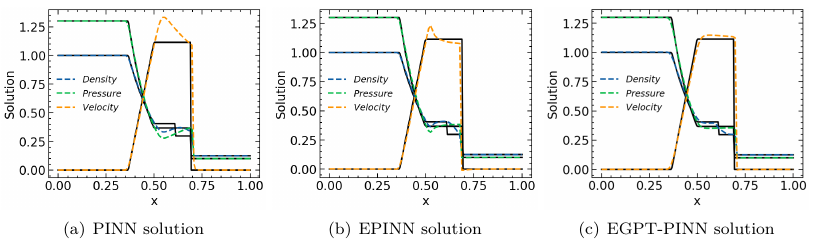}
    \caption{The predicted solutions for the inverse problem, obtained using the PINN, EPINN and EGPT-PINN methods, are presented in the second row of \Cref{tab:inverse_1}.}
    \label{fig:EGPT-inverse}
\end{figure}

\TBD{To test our inverse proglem solver, we select some $(x^*,t^*)$ and adopt the exact solution at the corresponding location as input data $\mathbf{U}^*$. 
We then train PINN (without RH-condition), EPINN and EGPT-PINN separately to predict the initial pressure $p_1$. The results are shown in \Cref{tab:inverse_1} and \Cref{fig:EGPT-inverse}. The EGPT-PINN demonstrates remarkable effectiveness in accurately inferring the unknown initial condition, as evidenced by the significantly superior performance of both the predicted $p_1$ and the simultaneously obtained solution compared to those produced by the PINN and EPINN.}

\begin{table}[htb]
\centering
\begin{tabular}{|c|c|c|c|c|}
\hline
$(x^*,t^*)$ & Reference $p_1$ & PINN & EPINN  &EGPT-PINN\\
\hline
(0.55,0.1)&1.3000&1.2995&1.2989&\textbf{1.2998}\\
\hline
(0.25,0.1)&1.3000&1.2993&1.2962&\textbf{1.2995}\\
\hline
(0.55,0.05)&1.6000&1.6003&1.5948&\textbf{1.6005}\\
\hline
(0.25,0.05)&1.6000&1.5990&1.5993&\textbf{1.5998}\\
\hline
\end{tabular}
\caption{Four test cases for the inverse 1D Sod problem by the PINN, EPINN and EGPT-PINN.}
\label{tab:inverse_1}
\end{table}

\section{Conclusion}\label{sec:conclusion}
In this paper, we develop the EGPT-PINN framework for the parameterized nonlinear conservation law problems and apply it to the inviscid Burgers' equation and the Euler equations. Our approach extends the capabilities of traditional PINNs and the recently introduced GPT-PINN by incorporating nonlinear model reduction techniques and entropy awareness while maintaining an unsupervised learning structure.

For Burgers' equation, the EGPT-PINN uses only one or two neurons to handle the discontinuities inherent in the shock waves, demonstrating superior accuracy and efficiency compared to traditional methods. In the case of the Euler equations, the EGPT-PINN's two-stage training process, with differentiated learning rates for the output and transform layers, ensured robust and rapid approximation of the solutions. This approach successfully managed the challenges posed by the shock waves and shock contact, with only 3 to 5 neurons, accurately capturing the dynamics of shock wave propagation and interaction. In both cases, the framework's ability to adaptively learn the system's parametric dependencies and incrementally expand the hidden layers proved crucial in capturing the complex behavior of shock interactions.

\bibliographystyle{abbrv}

\end{document}